\documentclass[12pt]{amsart}
\usepackage{righttag,epsf,amsfonts,amscd,latexsym}
\def\hepsffile{\leavevmode\epsffile}

\title[Gorenstein models...]
{Gorenstein models of del Pezzo surfaces of degree 1 over Dedekind schemes.}
\thanks{This work was partially supported by the grants RFBR no. 99--01--01132,
Grant of Leading Scientific Schools no. 96--15--96146, and INTAS-OPEN 97/2072}
\author{Mikhail Grinenko}
\address{Steklov Mathematical Institute}
\address{Max-Planck Institut f\"ur Mathematik}
\email{grin@mi.ras.ru / grinenko@mpim-bonn.mpg.de}
\date{}

\newtheorem{theorem}{\sc Theorem}[section]
\newtheorem{proposition}[theorem]{\sc Proposition}
\newtheorem{lemma}[theorem]{\sc Lemma}
\newtheorem{corollary}[theorem]{\sc Corollary}
\newtheorem{claim}[theorem]{\sc Claim}

\makeatletter

\@addtoreset{equation}{section}
\newcommand{\l@abcd}[2]{\hbox to\textwidth{#1\dotfill #2}}

\makeatother

\newcommand*{\mybegintheorem}[1]{\begin{trivlist}\it%
      \item[\hspace{\labelsep}{\bf #1}]}
\newcommand*{\myendtheorem}{\end{trivlist}}
\newenvironment*{theorem*}{\mybegintheorem{Theorem.}}{\myendtheorem}
\newenvironment*{proposition*}{\mybegintheorem{Proposition.}}{\myendtheorem}
\newenvironment*{corollary*}{\mybegintheorem{Corollary.}}{\myendtheorem}
\newenvironment*{definition*}{\mybegintheorem{Definition.}}{\myendtheorem}

\theoremstyle{remark}
\newtheorem{remark}[theorem]{\sc Remark}

\renewcommand{\phi}{\varphi}
\renewcommand{\epsilon}{\varepsilon}

\newcommand{\PT}{{{\mathbb P}^3}}
\newcommand{\PTw}{{{\mathbb P}^2}}

\newcommand{\PA}{{\mathbb P}}
\newcommand{\mc}{\mathcal}

\newcommand{\Bir}{\mathop{\rm Bir}\nolimits}

\newcommand{\PR}[1]{\mathbin{\mathop{\otimes}_{#1}}}

\newcommand{\Spec}{\mathop{{\bf Spec}\:}\nolimits}

\begin{document}
\abstract{Let $R$ be a Dedekind scheme, $\eta$ its generic point, $X$ and $V$
del Pezzo surfaces of degree 1 over $R$ that are Gorenstein Mori fiber spaces 
(as 3-folds germs over the ground field). We study birational maps
$\phi:X\dasharrow V$ over $R$ which are isomorphisms over the generic 
point of $R$. We put down normal forms of such transformations (in suitable 
coordinates) and give some properties of $X$ and $V$. In particular, we 
prove the uniqueness of a smooth model (corollary \ref{corol_smooth}).}
\endabstract
\maketitle

\section{Preliminary.}
\label{sec1}

The aim of this note is the following. Running the Sarkisov program
in order to factorize a birational map between two Mori fibrations by 
elementary links (see, for instance, \cite{BM}), we have to consider
birational maps that induce an isomorphism of fibers over generic points.
This is especially important in the birational rigidity problem.
It is amazing but during a long time we did not know nearly anything about
such transformations in cases of 3-fold fibrations on del Pezzo surfaces
of degree 1,2 or 3, except, probably, the Corti paper \cite{Corti}. Here, arguing in "coordinates" almost by the same way as in \cite{Corti},
we try to clear this subject in the practically important case of
Gorenstein fibrations on del Pezzo surfaces of degree 1. Our argumentation
is rather elementary and allows us to construct various interesting
examples.

Everywhere in the sequel the following conditions hold. The 
characteristic of the ground field $k$ is equal to 0. Let ${\mc O}$ be 
a DVR, $R=\Spec{\mc O}$, $K$ the field of functions of ${\mc O}$,
$X$ and $V$ del Pezzo surfaces of degree 1 over $R$. We assume $X$ and $V$
to be Gorenstein Mori fibrations over $R$. In particular, the fibers 
$X_K$ and $V_K$ over the generic point of $R$ have the Picard number 1.
We denote $X_0$ and $V_0$ the central fibers,
$$
    \phi: X\dasharrow V
$$
a birational map over $R$ inducing an isomorphism $\phi_K: X_K\simeq V_K$.

\begin{remark}
Recently in Park's preprint \cite{Park} it was proved that for degree of
fibers up to 4, $\phi$ is an isomorphism always if the central fibers 
are smooth.
\end{remark}

It is well known that $|-K_X|$ embeds $X$ into the weighted projective
space $\PA_{\mc O}(1,1,2,3)$ over ${\mc O}$, so that $X_0$ avoids the singular
points of $\PA_k(1,1,2,3)$ because of the Gorenstein property. As to $X_0$
itself, there are two possibilities (\cite{HiWa}, \cite{Reid1}):

\begin{proposition}
If $X_0$ is normal, then it has either a unique singular point which is
 minimally elliptic, or at most du Val singularities.

If $X$ is not normal, then its normalization is $\PTw$,  and $X$ is obtained
by gluing $\PTw$ along a quadric (possibly, reducible or non-reduced).
\end{proposition}

Note that in both the cases there exists a degree 2 morphism from $X$
onto a quadratic cone in $\PT$ branched along a cubic section that does not
pass through the cone vertex. In fact, this morphism is the restriction of 
a projection $\PA(1,1,2,3)\to \PA(1,1,2)$ to $X$ (the last weighted space
is nothing but a quadratic cone).

In order to separate terminal cases, we will need the following statement
(\cite{Reid2}):

\begin{proposition}
\label{th_Reid}
Let $U$ be a germ of a 3-fold terminal singularity, $S\in|-K_U|$ a general
element. Then $S$ has at most du Val singularities.
\end{proposition}

This statement will be used together with the so-called 
"recognition principle" (see \cite{BW}), which allows us to discern du Val singularities by an equation not in normal form.

The remaining part of the paper is organized as follows. First, we give
a suitable system of coordinates in $\PA(1,1,2,3)$, that allows us to
define equations for $\phi$, $X$, and $V$ in a convenient form (section
\ref{sec2}). In section \ref{sec3}, we produice the main division into
cases and deal with each of them. Then, we give some examples in 
section \ref{sec4}. Finally, section \ref{sec5} contains some remarks
that are closely related to the subject.

I would like to thank Prof. V.Shokurov and J.Park for useful discussions
of the problem. This work was carried out during
my stay at the Max-Planck Institute in Bonn. I would like also to express
my gratitude to the Directors and the stuff of the insitute for
hospitality and very nice atmosphere for a good work.

\section{Suitable coordinates.}
\label{sec2}

We denote ${\mathfrak m}$ the maximal ideal of ${\mc O}$. We may suppose
that it is generated by $t$, i.e., ${\mathfrak m}=(t){\mc O}$.

Then, we fix a copy of a weighted projective space 
$P_{1_{\mc O}}\cong\PA_{\mc O}(1,1,2,3)$ to which $X$ is embedded. 
By $P_1$ and $P_{1_K}$ we denote the specializations of 
$P_{1_{\mc O}}$ over the central point $t=0$ and the generic point 
respectively. Let $(x:y:z:w)$ be the coordinates in $P_{1_{\mc O}}$ of 
weights (1,1,2,3) (we will use the same denotations of coordinates for the
specializations). $P_{2_{\mc O}}$, $P_2$, $P_{2_K}$, and $(p:q:r:s)$
are used for $V$ in the same sense.

Further, $X$ is defined by a homogeneous polynomial of degree 6 in 
$P_{1_{\mc O}}$: 
$$
   aw^2+bz^3+cwzf_1(x,y)+z^2f_2(x,y)+zf_4(x,y)+f_6(x,y)=0,
$$
where $a,b,c\in{\mc O}$; $f_i$ denotes a homogeneous polynomial of degree $i$.
But $X_0$ does not pass through points $(0:0:1:0)$ and $(0:0:0:1)$, so
$a$ and $b$ are invertible in ${\mc O}$, and we suppose $a=1$. As to $b$, 
we may also assume $b=1$ since it is not important in the sequel, or,
if it is more convenient for the reader, we may extend ${\mc O}$ by the
cubic root of $b$ and put $z:=z\sqrt[3]{b}$. Note that such a substitution
can not increase the Picard number of $X_K$, since the new $K$ will be
an extension of degree 3, i.e., not even.

Now it is clear that substituting the coordinates $w$ and $z$ by a suitable
manner, we obtain the following equation for $X$:
\begin{equation}
\label{eq_X}
     w^2+z^3+zf_4(x,y)+f_6(x,y)=0
\end{equation}
($f_i$'s may be different from the previous ones, of course).

Further, since for $1\le i\le 3$ 
$H^1(P_{1_K},{\mc I}_{X_K}\times{\mc O}(i))=0$, 
we see that
$$
    H^0(X_K, -iK_{X_K})\simeq H^0(P_{1_K}, {\mc O}(i)).
$$
So the isomorphism $\phi^{-1}_K: V_K\simeq X_K$ induces
isomorphisms 
$$
   (\phi^{-1}_K)^*: H^0(P_{1_K}, {\mc O}(i))\simeq H^0(P_{2_K}, {\mc O}(i)),
$$
which yields an isomorphism
$$
    \phi^{-1}_K: P_{2_K}\simeq P_{1_K}.
$$

Thus, for the coordinates in $P_{2_K}$ we may put 
$$
\left\{
\begin{array}{ccl}
p & = & (\phi^{-1}_K)^*(x) \\
q & = & (\phi^{-1}_K)^*(y) \\
r & = & (\phi^{-1}_K)^*(z) \\
s & = & (\phi^{-1}_K)^*(w) 
\end{array}
\right.
$$

It remains to take into account a projection 
$$
      P_{2_K}=P_{2_{\mc O}}\PR{\mc O}\Spec K \to P_{2_{\mc O}},
$$ 
and then, multiplying by a suitable element of ${\mc O}$, 
for the coordinates in $P_{2_{\mc O}}$ we get $p=Ax$, $q=By$,
$r=Cz$, and $s=Dw$, where $A,B,C,D\in{\mc O}$. Now it is clear that we 
may suppose (possibly, stretching the coordinates by invertible elements)
\begin{equation}
\label{eq_phi}
\phi=
\left\{
\begin{array}{ccl}
p & = & t^ax \\
q & = & t^by \\
r & = & t^cz \\
s & = & t^dw 
\end{array}
\right\}
\end{equation}
where the set $(a,b,c,d)$ contains at least one 0. It is obvious that
\begin{equation}
\label{eq_phi-1}
\phi^{-1}=
\left\{
\begin{array}{ccl}
x & = & t^{\alpha}p \\
y & = & t^{\beta}q \\
z & = & t^{\gamma}r \\
w & = & t^{\delta}s 
\end{array}
\right\}
\end{equation}
with the same condition on the set $(\alpha,\beta,\gamma,\delta)$. Moreover,
taking into account the weights, we get
$$
   6(a+\alpha)=6(b+\beta)=3(c+\gamma)=2(d+\delta).
$$

Then, substituting (\ref{eq_phi-1}) for the coordinates in (\ref{eq_X})
and using the condition that $V_0$ does not pass through points (0:0:1:0)
and (0:0:0:1), we obtain an equation of $V$
\begin{equation}
\label{eq_V}
     s^2+r^3+rg_2(p,q)+g_6(p,q)=0
\end{equation}
and a condition $3\gamma=2\delta$, where $g_i$ are some homogeneous 
polynomials of degree $i$.

Finally, we see that there exists a positive integer $m$ such that
the following conditions hold:
\begin{equation}
\label{conds}
\left\{
\begin{array}{rcl}
a+\alpha & = & m  \\
b+\beta  & = & m  \\
c+\gamma & = & 2m \\
d+\delta & = & 3m \\
2d       & = & 3c \\
2\delta  & = & 3\gamma
\end{array}
\right.
\end{equation}

\section{Main division into cases.}
\label{sec3}

Using the results of section \ref{sec2} and the symmetry of the situation, 
we can produice the following division into four cases:
$$
\begin{array}{|c|c|c|c|}
\hline
\mbox{Case} & (a,b,c,d) & (\alpha,\beta,\gamma,\delta) & \mbox{Remarks} \\
\hline
A & (a,m,0,0) & (\alpha,m,2m,3m) & \begin{array}{c}
                                   a+\alpha=m \\
                                   a,\alpha>0
                                 \end{array}       \\
\hline
B & (0,m,0,0) & (m,0,2m,3m) &  \\
\hline
C & (m,m,0,0) & (0,0,2m,3m) &  \\
\hline
D & (0,m,2k,3k) & (m,0,2l,3l) & \begin{array}{c}
                                 k+l=m   \\
                                 m\ge 2  \\
                                 k,l>0  
                                \end{array}       \\
\hline
\end{array}
$$

We will show that 
\begin{center}
{\it The cases A, B, and C can not occur.}
\end{center}

Here our main tools will be proposition \ref{th_Reid}, the "recognition
principle", and a condition that the central fibers are Gorenstein.

First, assuming $X$ and $V$ to be defined by equations (\ref{eq_X}) and
(\ref{eq_V}) while $\phi$ and $\phi^{-1}$ are in form (\ref{eq_phi}) and
(\ref{eq_phi-1}), we obtain the following conditions:

\begin{equation}
\label{cond_eq}
\begin{array}{c}
\left\{
 \begin{array}{ccl}
  f_4(x,y) & = & t^{c-2d}g_4(t^ax,t^by) \\
  f_6(x,y) & = & t^{-2d}g_6(t^ax,t^by)
 \end{array}
\right.  \\
 \\
\left\{
 \begin{array}{ccl}
  g_4(p,q) & = & t^{\gamma-2\delta}f_4(t^{\alpha}x,t^{\beta}y) \\
  g_6(p,q) & = & t^{-2\delta}f_6(t^{\alpha}x,t^{\beta}y)
 \end{array}
\right.  \\
\end{array}
\end{equation}

Then, we will suppose that $g_4$ and $g_6$ are defined by
\begin{equation}
\label{eq_b}
\begin{array}{rclcrcl}
g_4(p,q) & = & \sum_{i=0}^4a_ip^{4-i}q^i; &&
g_6(p,q) & = & \sum_{i=0}^6b_ip^{6-i}q^i. \\
\end{array}
\end{equation}

\noindent{\bf Case A.} Here we have
$$
\begin{array}{rcccl}
f_4(x,y) & = & \sum_{i=0}^4a_it^{4a+(m-a)i}x^{4-i}y^i & = & 
t^{4a}g_4(x,t^{m-a}y), \\
f_6(x,y) & = & \sum_{i=0}^6b_it^{6a+(m-a)i}x^{6-i}y^i & = & 
t^{6a}g_6(x,t^{m-a}y). \\
\end{array}
$$
So $X$ is defined by
$$
w^2+z^3+t^{4a}zg_4(x,t^{m-a}y)+t^{6a}g_6(x,t^{m-a}y)=0.
$$

Since $a>0$, the curve $\{t=w=z=0\}$ lies on $X$, and $X$ is singular along it.
So it can not be a terminal case.

\noindent{\bf Case C.} The same reason as above: simply substitute $m$ for $a$.

\noindent{\bf Case B.} We have
$$
\begin{array}{rcl}
f_4(x,y) & = & \sum_{i=0}^4a_it^{mi}x^{4-i}y^i, \\
f_6(x,y) & = & \sum_{i=0}^6b_it^{mi}x^{6-i}y^i,\\
\end{array}
$$
and $X$ is
$$
w^2+z^3+z(\sum_{i=0}^4a_it^{mi}x^{4-i}y^i)+\sum_{i=0}^6b_it^{mi}x^{6-i}y^i=0.
$$

Note that at least one of the coefficients is not identically 0.
Choose an affine piece $\{y\ne 0\}$. Using the same denotations for the
coordinates in this affine piece, we see that a point $A=\{t=x=z=w=0\}$ is
singular on $X$.

Then, a general element of $|-K_X|$ through $A$ has the form $\{x=th\}$ in
the affine part, where $h\in{\mc O}$. Thus, we obtain that such an element
is defined by
$$
w^2+z^3+uzt^4+vt^6=0,
$$
where $u,v\in{\mc O}$ and at least one of them is not 0. Due the "recognition
principle", it can not be an equation of a du Val singularity. By proposition
\ref{th_Reid}, $X$ can not be terminal.

So we proved that the cases A,B, and C can not occur. Now we shall deal with
the case D.

\medskip

\noindent{\bf The unique possibility: case D.} In the sequel we will always
assume that
$$
\begin{array}{cc}
\phi=
\left\{
\begin{array}{ccl}
p & = & x \\
q & = & t^my \\
r & = & t^{2k}z \\
s & = & t^{3k}w 
\end{array}
\right\},
&
\phi^{-1}=
\left\{
\begin{array}{ccl}
x & = & t^mp \\
y & = & q \\
z & = & t^{2l}r \\
w & = & t^{3l}s 
\end{array}
\right\},
\end{array}
$$
where $k+l=m$, and $0<k\le l$, so $k\le\frac{m}2$.

Using (\ref{cond_eq}), we see that
$$
\begin{array}{lcr}
  f_4(x,y)=t^{-4k}g_4(x,t^my), & & f_6(x,y)=t^{-6k}g_6(x,t^my).
\end{array}
$$
Thus, there exist $\alpha_0,\alpha_1,\beta_0,\beta_1,\beta_2\in{\mc O}$
such that
\begin{equation}
\label{cond_coeff}
\left\{
\begin{array}{rcl}
a_0 & = & \alpha_0t^{4k} \\
a_1 & = & \alpha_1t^{4k-m} \\
b_0 & = & \beta_0t^{6k} \\
b_1 & = & \beta_1t^{6k-m} \\
b_2 & = & \beta_2t^{6k-2m} \\
\end{array}
\right.
\end{equation}

For convenience, let us put down all equations that we will need in the
sequel:
\begin{equation}
\label{eq_fg}
\begin{array}{rcl}
f_4(x,y) & = & \alpha_0x^4+\alpha_1x^3y+a_2t^{2m-4k}x^2y^2+
          a_3t^{3m-4k}xy^3+ \\
&&        +a_4t^{4m-4k}y^4, \\
&& \\
f_6(x,y) & = & \beta_0x^6+\beta_1x^5y+\beta_2x^4y^2+
          b_3t^{3m-6k}x^3y^3+ \\
&&   +b_4t^{4m-6k}x^2y^4+b_5t^{5m-6k}xy^5+b_6t^{6m-6k}y^6, \\
&& \\
g_4(p,q) & = & \alpha_0t^{4k}p^4+\alpha_1t^{4k-m}p^3q+a_2p^2q^2+a_3pq^3+
   a_4q^4, \\
&& \\
g_6(p,q) & = & \beta_0t^{6k}p^6+\beta_1t^{6k-m}p^5q+\beta_2t^{6k-2m}p^4q^2+
   b_3p^3q^3+ \\
&&  b_4p^2q^4+b_5pq^5+b_6q^6, \\
\end{array}
\end{equation}
and
\begin{equation}
\label{eq_XV}
\begin{array}{cl}
X: &  \left\{
 \begin{array}{l}
   w^2+z^3+z\left(\alpha_0x^4+\alpha_1x^3y+
    t^{2m-4k}\sum_{i=2}^4a_it^{m(i-2)}x^{4-i}y^i\right)+ \\
 \\
   +\beta_0x^6+\beta_1x^5y+\beta_2x^4y^2+
      t^{3m-6k}\sum_{i=3}^6b_it^{m(i-3)}x^{6-i}y^i=0, \\
 \end{array} 
\right. \\
& \\
V: &  \left\{
 \begin{array}{l}
  s^2+r^3+r\left(\alpha_0t^{4k}p^4+\alpha_1t^{4k-m}p^3q+
  \sum_{i=2}^4a_ip^{4-i}q^i\right)+ \\
 \\
  +\beta_0t^{6k}p^6+\beta_1t^{6k-m}p^5q+\beta_2t^{6k-2m}p^4q^2+
     \sum_{i=3}^6b_ip^{6-i}q^i=0.
 \end{array} 
\right. \\
\end{array}
\end{equation}

Note that $\phi$ is not defined along the curve $\{t=x=0\}$,
and $\phi^{-1}$ along $\{t=q=0\}$. Then, $V_0$ is contracted to
the point $A=\{t=0,(0:1:0:0)\}\in X$, and $X_0$ to the point
$B=\{t=0,(1:0:0:0)\}\in V$.

\begin{claim}
$X$ is always singular.
\end{claim}
\noindent{\sc Proof.} Since $2m-4k\ge 0$, it is easy to check that the 
point $A$ is singular.

In particular, we proved the following important result:

\begin{corollary}
\label{corol_smooth}
Let $U/T$ be a Mori fibration on del Pezzo surfaces of degree 1 over 
a curve $T$. Denote
$$
  {\mc {MF}}_{base}=\left\{ \mbox{Mori fibrations that are birational over $T$
to $U$}\right\}.
$$
Then ${\mc {MF}}_{base}$ contains at most one non-singular element.
\end{corollary}

\begin{lemma}
\label{lem_m6k}
Always $m\le 6k$.
\end{lemma}
\noindent{\sc Proof.} Assume the converse, i.e., $m>6k$. Then 
$\alpha_1,\beta_1\in{\mathfrak m}$ and $\beta_2\in{\mathfrak m}^6$
because of (\ref{cond_coeff}). Let $S=\{x=th\}$, $h\in{\mc O}$, be a
general element of $|-K_X-A|$ (in the affine piece $y\ne 0$). Then
$S$ has the form
$$
w^2+z^3+uzt^4+vt^6=0
$$
for some $u,v\in{\mc O}$. So $S$ is not canonical, which contradicts to
proposition \ref{th_Reid}.

\medskip

Now we consider two possible cases: $k=l$, i.e., $m=2k$, and $k<l$.

\noindent{\bf Case ${\bf k=l}$.} In this case, let us consider the 
affine part $\{p\ne 0\}$ of $V$. Then $V$ is defined by
$$
s^2+r^3+r(\alpha_0t^{4k}+q(\ldots\mbox{terms}\ldots))+
 \beta_0t^{6k}+\beta_1t^{4k}q+q^2(\ldots\mbox{terms}\ldots)=0,
$$
and we see that the point $B\in V$ is singular. Thus, for $k=l$, 
such a fiber-wise transformations only exist between singular varieties.

\noindent{\bf Case ${\bf k<l}$.} Here, it would be interesting to know
under which conditions $V$ is non-singular. Since $X_0$ is contracted to
the point $B$, we should first check $V$ at this point. Take an affine 
piece $\{p\ne 0\}$. Then from equations (\ref{eq_XV}) it follows that
$\beta_1t^{6k-m}\not\in{\mathfrak m}$, i.e., $6k-m\ge 0$. Thus, $m=6k$
because of lemma \ref{lem_m6k}, and $\beta_1$ is invertible in ${\mc O}$.

So, if $V$ is smooth, then $(a,b,c,d)=(0,6k,2k,3k)$ and 
$(\alpha,\beta,\gamma,\delta)=(6k,0,10k,15k)$. Let us note that 
$\alpha_1\in{\mathfrak m}^{2k}$ and $\beta_2\in{\mathfrak m}^{6k}$ in
this case. It easy to check also that $X_0$ only has a $E_8$-singularity at
$A$, so $X$ has $cE_8$ there.

\section{Examples.}
\label{sec4}

Here we give examples with some interesting properties.

\noindent{\bf Example 1. Smooth case.} Suppose $(a,b,c,d)=(0,6,2,3)$ and
$(\alpha,\beta,\gamma,\delta)=(6,0,10,15)$. $X$ and $V$ is defined by
$$
\begin{array}{ccl}
X: && w^2+z^3+x^5y+t^{24}xy^5 = 0, \\
&&\\
V: && s^2+r^3+p^5q+pq^5 = 0.
\end{array}
$$
It is easy to see that $V$ is non-singular, $X$ has a $cE_8$-singularity
at the point $\{t=0, (0:1:0:0)\}$.

\medskip

\noindent{\bf Example 2. Birational automorphism.} Let $(a,b,c,d)=(2,0,2,3)$
and $(\alpha,\beta,\gamma,\delta)=(0,2,2,3)$. Equations:
$$
\begin{array}{ccl}
X: && w^2+z^3+t^4x^5y+xy^5 = 0, \\
&&\\
V: && s^2+r^3+p^5q+t^4pq^5 = 0.
\end{array}
$$
Note that $X$ and $V$ are isomorphic to each other: simply take $w=s$,
$z=r$, $x=q$, and $y=p$. So we may assume that $\phi\in\Bir(X)$ is defined
by a transformation of coordinates in $P_{1_K}$: 
$$
\begin{array} {rcl}
x & \longmapsto & t^{-1}y; \\
y & \longmapsto & tx; \\
z & \longmapsto & z; \\
w & \longmapsto & w.
\end{array}
$$
$X$ has a $cE_8$-singularity in the central fiber. 

\medskip

\noindent{\bf Example 3. Non-normal fibers.} We give two examples with
the same weights $(a,b,c,d)=(2,0,2,3)$ and
$(\alpha,\beta,\gamma,\delta)=(0,2,2,3)$ such that $X_0$ and $V_0$ are not
normal. Both the cases are also birational automorphisms. The first one is
$$
\begin{array}{ccl}
X: && w^2+z^3+tzx^2y^2+tx^6+t^7y^6 = 0, \\
&&\\
V: && s^2+r^3+trp^2q^2+t^7p^6+tq^6 = 0.
\end{array}
$$
$X_0$ is non-normal with an equation $w^2+z^3=0$ (the "cusp" case, see
\cite{Reid1}, 1.4.). $X$ has a $cE_7$-singularity in the central fiber.

The second is
$$
\begin{array}{ccl}
X: && w^2+z^3-3zx^2y^2+tx^6+2x^3y^3+t^7y^6 = 0, \\
&&\\
V: && s^2+r^3-3rp^2q^2+t^7p^6+2p^3q^3+tq^6 = 0.
\end{array}
$$
$X_0$ is defined by $w^2+(z-xy)^2(z+2xy)=0$ (the "node" case), and
$X$ is $cD_4$.

\section{Remarks.}
\label{sec5}

As it was shown, both $\phi$ and $\phi^{-1}$ contract central fibers
to points. It gives us some special conditions on central fibers and
elements of $|-K|$.

\medskip

\noindent{\bf Central fibers.} Since $X$ and $V$ are assumed to be 
terminal, i.e., having only rational singularities, their central 
fibers $X_0$ and $V_0$ must be rational. Consider $X_0$. If it is
not normal, then its normalization is $\PTw$, so $X_0$ is rational
(example 3 of the previous section). If $X_0$ is normal, then it is 
rational if and only if it has at most du Val singularities. So the
remaining case of minimally elliptic singularities (defined by 
$w^2+z(z^2-x^4)=0$) is not possible.

\medskip

\noindent{\bf Anticanonical divisors.} Let $S\in|-K_X|$ be a general 
element, $S'\in|-K_V|$ its strict transform. Then $S$ is a smooth
elliptic surface. We will show that the elliptic surface $S'$ has 
an elliptic (non-minimal in general case) singularity.

Denote $s$ and $s'$ sections of $S$ and $S'$ that are the base locus
of anticanonical divisors, $\rho:T\to S'$ the minimal resolution 
of singularities of $S'$. Since $S$ is relatively minimal over $R$,
it yields a birational morphism $\mu: T\to S$, which is factorized
by contractions of {\bf -1}-curves. Denote $C$ and $C'$ the central
fibers of $S$ and $S'$, $\tilde C$ and $\tilde C'$ their strict
transforms on $T$, and $s_T$ the pre-image of $s$ (or $s'$) on $T$. 
Obviously, $C$ is either an elliptic curve, or
a rational curve with a double point ("cusp" or "node"). We also know
that $\rho$ contracts $\tilde C$ to a point on $S'$, and $\mu$ contracts
$\tilde C'$. Then, since $C'$ is irreducible, $\tilde C'$ is a unique 
{\bf -1}-curve on $T$. Moreover, $s'\cap C'$ is a non-singular point on 
$S'$, thus $s_T$ does not intersect the exeptional curves. Further, the
point $s\cap C$ is a non-singular point of $C$, so $\tilde C$ has the
arithmetic genus 1. The following picture clarifies the geometric 
situation.

\begin{figure}[htbp]
 \begin{center}
  \epsfxsize 10cm
  \hepsffile{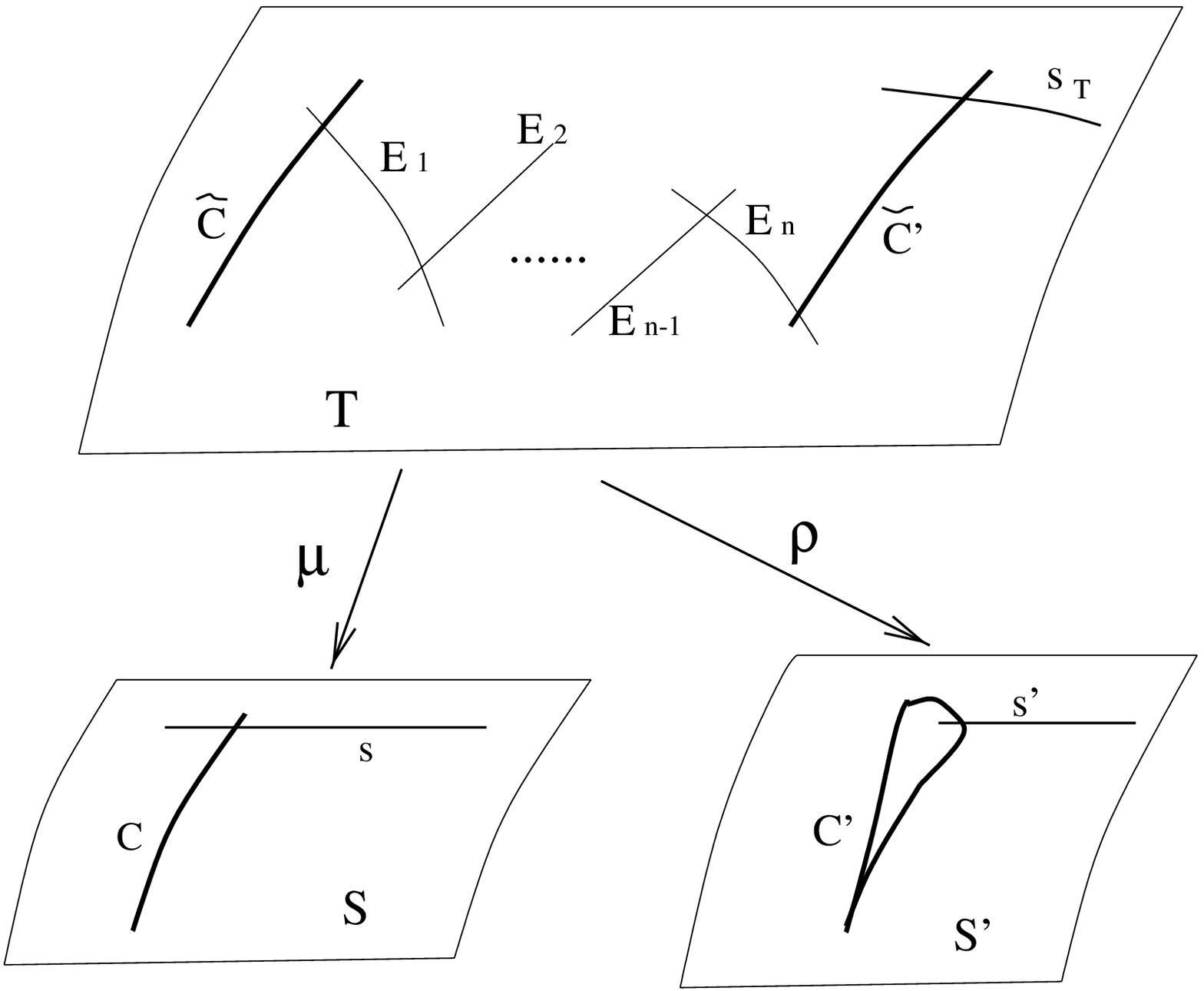} 
 \end{center}
\caption{}\label{kuku.fig}
\end{figure}

Here $E_1,\ldots,E_n$, $n\ge 0$, are {\bf -2}-curves. Let us note also 
that $C'$ can not have a node point, so it is a cusp.

Thus, we have the following intersection numbers: 
$\tilde C^2=\tilde C'{}^2=-1$, $E_i^2=-2$, $E_i\circ E_{i+1}=1$. 

Then, comparing $K_T$ for morphisms $\rho$ and $\mu$, we see that
$$
K_T\sim -(n+1)\tilde C-nE_1-(n-1)E_2-\ldots-E_n
$$
and
$$
    K_T\sim E_1+2E_2+\ldots nE_n+(n+1)\tilde C',
$$
whence
$$
(n+1)(\tilde C+E_1+\ldots+E_n+\tilde C')\sim 0
$$
(in fact, it is obvious, since $\tilde C+E_1+\ldots+E_n+\tilde C'$ is a 
fiber of $T$). This formula can be explained as follows. Suppose for 
instance that $S$ and $S'$ are projective elliptic surfaces, $f$ and $f'$
are their classes of a fiber. There exist integers $m$ and $m'$ such
that $K_S\sim mf$ and $K_{S'}\sim m'f'$. Then $m'=m+(n+1)$.

Now it is clear that $S'$ has an elliptic singularity with $\tilde C$ as
the fundamental cycle, so it is minimally elliptic if and only if $n=0$.


\begin{thebibliography}{10}

\bibitem{BM}
Bruno, A., and Matsuki, K. {\it Log Sarkisov program}. 
Int. J. Math. {\bf 8} (1997), No.4, 451-494 

\bibitem{BW}
Bruce, J.W., and Wall, C.T.C. {\it On the classification of cubic 
surfaces}. J. London Math. Soc., II. Ser. 19 (1979), 245-256

\bibitem{Corti}
Corti, A., {\it Del Pezzo surfaces over Dedekind schemes}. 
Ann. Math., II. Ser. 144, {\bf 3} (1996), 641-683

\bibitem{HiWa}
Hidaka, F., and Watanabe, K.-i. {\it Normal Gorenstein surfaces with ample
anti-canonical divisor}. Tokyo J. Math., {\bf 4} (1981), No. 2, 319-330


\bibitem{Park}
Park, J., {\it Birational maps of del Pezzo fibrations}. E-preprint
math/9912076


\bibitem{Reid1}
Reid, M., {\it Nonnormal del Pezzo surfaces}.
Publ. Res. Inst. Math. Sci. {\bf 30} (1994), No.5, 695-727

\bibitem{Reid2}
Reid, M., {\it Young person's guide to canonical singularities}.
Algebraic geometry, Proc. Summer Res. Inst., Brunswick/Maine 1985, 
part 1, Proc. Symp. Pure Math. {\bf 46} (1987), 345-414


\end{thebibliography}
\end{document}